\documentstyle[12pt, leqno, amstex]{amsart}

\makeatletter

\begin{document}
\bibliographystyle{amsalpha}

\newtheorem{Assumption}{Assumption}[section]
\newtheorem{Theorem}{Theorem}[section]
\newtheorem{Lemma}{Lemma}[section]
\newtheorem{Remark}{Remark}[section]
\newtheorem{Corollary}{Corollary}[section]
\newtheorem{Conjecture}{Conjecture}[section]
\newtheorem{Proposition}{Proposition}[section]
\newtheorem{Example}{Example}[section]
\newtheorem{Definition}{Definition}[section]
\newtheorem{Question}{Question}[section]
\renewcommand{\thesubsection}{\it}

\baselineskip=14pt
\addtocounter{section}{-1}

\title{Symmetries in the system of type $A_5^{(2)}$}

\author{Yusuke Sasano}

\address{Graduate School of Mathematical 
Sciences, The University of Tokyo, Tokyo 153-8914, Japan.}
\email{sasano@ms.u-tokyo.ac.jp}

\begin{abstract}
In this paper, we propose a 3-parameter family of coupled Painlev\'e III systems in dimension four with affine Weyl group symmetry of type $A_5^{(2)}$. We also propose its symmetric form in which the $A_5^{(2)}$-symmetries become clearly visible.
\end{abstract}
\maketitle

\section{Statement of main results} 
In \cite{Sasa5,Sasa6,Sasa8}, we presented some types of coupled Painlev\'e systems with various affine Weyl group symmetries. In this paper, we present a 3-parameter family of coupled Painlev\'e III systems with affine Weyl group symmetry of type $A_5^{(2)}$. This system is the first example with the B{\"a}cklund transformations satisfying Noumi-Yamada's universal description for $A_5^{(2)}$ root system (see \cite{N2}). At first, we propose a 3-parameter family of autonomous ordinary differential systems with the invariant divisors $f_i$ as variables:
\begin{equation}
  \left\{
  \begin{aligned}
   \frac{df_0}{dt} &=\{2(f_1-f_3)g_1-2g_2f_1-\alpha_1-\alpha_3\}f_0+\alpha_0(f_1-f_3),\\
   \frac{df_1}{dt} &=\{2(f_0-f_3)g_2-2g_1f_0-\alpha_0-\alpha_3\}f_1+\alpha_1(f_0-f_3),\\
   \frac{df_2}{dt} &=\{(f_3-f_1+3f_0)g_1+(f_3+3f_1-f_0)g_2+1\}f_2-4\alpha_2 g_1g_2,\\
   \frac{df_3}{dt} &=-(2f_0g_1+2f_1g_2+\alpha_0+\alpha_1)f_3-\alpha_3(f_0+f_1),\\
   \frac{dg_1}{dt} &=(f_0-f_1+f_3)g_1^2+\{(f_0-f_1-f_3)g_2+\alpha_0+\alpha_1+\alpha_3\}g_1+f_2(f_3+3f_1-f_0),\\
   \frac{dg_2}{dt} &=(f_1-f_0+f_3)g_2^2+\{(f_1-f_0-f_3)g_1+\alpha_0+\alpha_1+\alpha_3\}g_2+f_2(f_3+3f_0-f_1).\\
   \end{aligned}
  \right. 
\end{equation}
Here $f_0,f_1,f_2,f_3$ and $g_1,g_2$ denote unknown complex variables and $\alpha_0,\dots,\alpha_3$ are the parameters satisfying the condition
$$
\alpha_0+\alpha_1+2\alpha_2+\alpha_3=\frac{1}{2}.
$$

\begin{Proposition}
This system has the following invariant divisors\rm{:\rm}
\begin{center}
\begin{tabular}{|c|c|c|} \hline
invariant divisors & parameter's relation \\ \hline
$f_0:=0$ & $\alpha_0=0$  \\ \hline
$f_1:=0$ & $\alpha_1=0$  \\ \hline
$f_2:=0$ & $\alpha_2=0$  \\ \hline
$f_3:=0$ & $\alpha_3=0$  \\ \hline
\end{tabular}
\end{center}
\end{Proposition}

\begin{Theorem}
This system is invariant under the transformations $s_0,s_1,s_2,s_3,\pi$ defined as follows$:$ with {\it the notation} $(*):=(f_0,f_1,f_2,f_3,g_1,g_2;\alpha_0,\alpha_1,\alpha_2,\alpha_3),$
\begin{align}
\begin{split}
s_0:(*) \rightarrow &(f_0,f_1,f_2+\frac{\alpha_0g_2}{f_0},f_3,g_1+\frac{\alpha_0}{f_0},g_2;-\alpha_0,\alpha_1,\alpha_2+\alpha_0,\alpha_3),\\
s_1:(*) \rightarrow &(f_0,f_1,f_2+\frac{\alpha_1g_1}{f_1},f_3,g_1,g_2+\frac{\alpha_1}{f_1};\alpha_0,-\alpha_1,\alpha_2+\alpha_1,\alpha_3),\\
s_2:(*) \rightarrow &(f_0-\frac{\alpha_2g_2}{f_2},f_1-\frac{\alpha_2g_1}{f_2},f_2,f_3-\frac{\alpha_2(g_1+g_2)}{f_2},g_1,g_2;\\
&\alpha_0+\alpha_2,\alpha_1+\alpha_2,-\alpha_2,\alpha_3+2\alpha_2),\\
s_3:(*) \rightarrow &(f_0,f_1,f_2+\frac{\alpha_3(g_1+g_2)}{f_3}+\frac{\alpha_3^2}{f_3^2},f_3,g_1+\frac{\alpha_3}{f_3},g_2+\frac{\alpha_3}{f_3};\\
&\alpha_0,\alpha_1,\alpha_2+\alpha_3,-\alpha_3),\\
\pi:(*) \rightarrow &(f_1,f_0,f_2,f_3,g_2,g_1;\alpha_1,\alpha_0,\alpha_2,\alpha_3).
\end{split}
\end{align}
\end{Theorem}
Here the Poisson bracket $\{,\}$ is defined by
$$
\{f_2,f_3\}=g_1+g_2, \ \{f_3,g_1\}=\{f_3,g_2\}=1.
$$

\begin{Theorem}
This system has two first integrals\rm{:\rm}
$$
\frac{df_3}{dt}=\frac{d(f_0+f_1)}{dt}, \quad \frac{d(f_2-g_1g_2)}{dt}=f_2-g_1g_2.
$$
\end{Theorem}
From this, we have
$$
f_3=f_0+f_1-1, \quad f_2-g_1g_2=e^{(t+c)}.
$$
Here we set
$$
t+c=log T, \quad q_1:=g_1, \ p_1:=f_0, \ q_2:=g_2, \ p_2:=f_1,
$$
then we obtain a 3-parameter family of coupled Painlev\'e III systems in dimension four with affine Weyl group symmetry of type $A_5^{(2)}$ explicitly given by
\begin{equation}
  \left\{
  \begin{aligned}
   \frac{dq_1}{dT} &=\frac{2q_1^2p_1-q_1^2+(\alpha_0+\alpha_1+\alpha_3)q_1}{T}-1+4p_2+\frac{2q_1q_2p_2}{T},\\
   \frac{dp_1}{dT} &=\frac{-2q_1p_1^2+2q_1p_1-(\alpha_0+\alpha_1+\alpha_3)p_1+\alpha_0}{T}-\frac{2p_1q_2p_2}{T},\\
   \frac{dq_2}{dT} &=\frac{2q_2^2p_2-q_2^2+(\alpha_0+\alpha_1+\alpha_3)q_2}{T}-1+4p_1+\frac{2q_1p_1q_2}{T},\\
   \frac{dp_2}{dT} &=\frac{-2q_2p_2^2+2q_2p_2-(\alpha_0+\alpha_1+\alpha_3)p_2+\alpha_1}{T}-\frac{2q_1p_1p_2}{T}\\
   \end{aligned}
  \right. 
\end{equation}
with the Hamiltonian
\begin{align}
\begin{split}
H&=\frac{q_1^2p_1^2-q_1^2p_1+(\alpha_0+\alpha_1+\alpha_3)q_1p_1-\alpha_0q_1}{T}-p_1\\
&+\frac{q_2^2p_2^2-q_2^2p_2+(\alpha_0+\alpha_1+\alpha_3)q_2p_2-\alpha_1q_2}{T}-p_2+4p_1p_2+\frac{2q_1p_1q_2p_2}{T}.
\end{split}
\end{align}
Here the Poisson bracket $\{,\}$ is defined by
$$
\{q_1,p_1\}=\{q_2,p_2\}=1, \ \{q_1,q_2\}=\{q_1,p_2\}=\{p_1,q_2\}=\{p_1,p_2\}=0.
$$

\begin{figure}
\unitlength 0.1in
\begin{picture}(30.20,21.70)(17.50,-24.70)
%
\special{pn 20}%
\special{ar 2100 650 350 350  0.0000000 6.2831853}%
%
\special{pn 20}%
\special{ar 2130 2120 350 350  0.0000000 6.2831853}%
%
\special{pn 20}%
\special{ar 3120 1380 350 350  0.0000000 6.2831853}%
%
\special{pn 20}%
\special{ar 4420 1390 350 350  0.0000000 6.2831853}%
%
\special{pn 20}%
\special{pa 2390 830}%
\special{pa 2850 1140}%
\special{fp}%
%
\special{pn 20}%
\special{pa 2380 1880}%
\special{pa 2850 1630}%
\special{fp}%
\put(20.1000,-22.4000){\makebox(0,0)[lb]{$p_1$}}%
\put(19.7000,-7.6000){\makebox(0,0)[lb]{$p_2$}}%
\put(28.7000,-14.9000){\makebox(0,0)[lb]{$q_1q_2+T$}}%
\put(40.7000,-14.5000){\makebox(0,0)[lb]{$p_1+p_2-1$}}%
%
\special{pn 20}%
\special{pa 4090 1230}%
\special{pa 3460 1230}%
\special{fp}%
\special{sh 1}%
\special{pa 3460 1230}%
\special{pa 3527 1250}%
\special{pa 3513 1230}%
\special{pa 3527 1210}%
\special{pa 3460 1230}%
\special{fp}%
%
\special{pn 20}%
\special{pa 4100 1560}%
\special{pa 3450 1560}%
\special{fp}%
\special{sh 1}%
\special{pa 3450 1560}%
\special{pa 3517 1580}%
\special{pa 3503 1560}%
\special{pa 3517 1540}%
\special{pa 3450 1560}%
\special{fp}%
\end{picture}%
\label{B3}
\caption{The transformations $s_i$ satisfy the relations: $s_i^2=1 \ (i=0,1,2,3), \ (s_0s_1)^2=(s_0s_3)^2=(s_1s_3)^2=1, \ (s_0s_2)^3=(s_1s_2)^3=1, \ (s_2s_3)^4=1.$}
\end{figure}
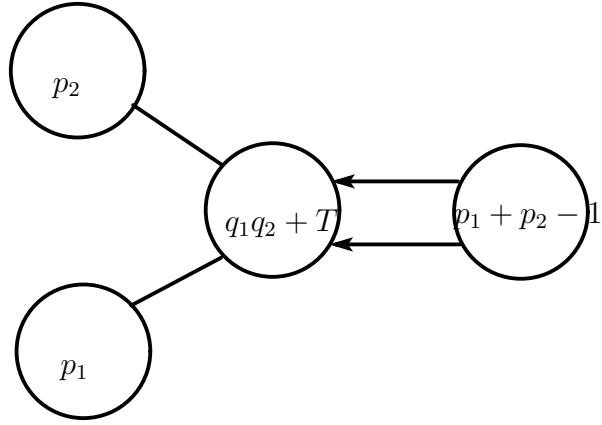

\begin{Theorem}
This system is invariant under the transformations $s_0,s_1,s_2,s_3,\pi$ defined as follows$:$ with {\it the notation} $(*):=(q_1,p_1,q_2,p_2,T;\alpha_0,\alpha_1,\alpha_2,\alpha_3),$
\begin{align}
\begin{split}
s_0:(*) \rightarrow &(q_1+\frac{\alpha_0}{p_1},p_1,q_2,p_2,T;-\alpha_0,\alpha_1,\alpha_2+\alpha_0,\alpha_3),\\
s_1:(*) \rightarrow &(q_1,p_1,q_2+\frac{\alpha_1}{p_2},p_2,T;\alpha_0,-\alpha_1,\alpha_2+\alpha_1,\alpha_3),\\
s_2:(*) \rightarrow &(q_1,p_1-\frac{\alpha_2 q_2}{q_1q_2+T},q_2,p_2-\frac{\alpha_2 q_1}{q_1q_2+T},T;\alpha_0+\alpha_2,\alpha_1+\alpha_2,-\alpha_2,\alpha_3+2\alpha_2),\\
s_3:(*) \rightarrow &(q_1+\frac{\alpha_3}{p_1+p_2-1},p_1,q_2+\frac{\alpha_3}{p_1+p_2-1},p_2,T;\alpha_0,\alpha_1,\alpha_2+\alpha_3,-\alpha_3),\\
\pi:(*) \rightarrow &(q_2,p_2,q_1,p_1,T;\alpha_1,\alpha_0,\alpha_2,\alpha_3).
\end{split}
\end{align}
\end{Theorem}

\end{document}